# Extensions and Applications of Stein–Weiss Operators to the Study of Traceless Symmetric Tensors

Sergey Stepanov[1, 2], Irina Tsyganok[2]

[1] Department of Scientific Information on Fundamental and Applied Mathematics, Russian Institute for Scientific and Technical Information of the Russian Academy of Sciences, 125190 Moscow, Russia,
E-mail: s.e.stepanov@mail.ru

[2] Dept. of Mathematics and Data Analysis, Finance University, 125468 Moscow, Russia,
E-mail: i.i.tsyganok@mail.ru

**Abstract**

First-order differential operators arising from the representation-theoretic decomposition of the covariant derivative play a fundamental role in Riemannian geometry. In this paper, we investigate Stein–Weiss $O(n)$-gradients acting on covariant symmetric trace-free tensors of arbitrary rank $p \geq 2$. By analyzing the decomposition of $T^*M \otimes S_0^p(M)$ into irreducible $O(n)$-modules, we explicitly describe the associated generalized gradients and derive Weitzenböck formulas for the compositions of these operators with their formal adjoints. Our results extend Bouguignon's four-dimensional formulas for the case $p = 2$ and generalize earlier work by other authors to symmetric tensors of higher rank. The resulting formulas provide a unified framework for the study of second-order Stein–Weiss operators and furnish analytical tools applicable to deformation complexes, curvature estimates, and stability problems in geometric analysis. This article continues the authors' previous investigations of Stein–Weiss operators on natural tensor bundles.

## 1. Introduction

First-order differential operators arising from the decomposition of the covariant derivative play a fundamental role in the analysis of geometric structures on smooth manifolds. Such decompositions of the covariant derivative into first-order differential operators are often associated with the structure group of the principal bundle of the manifold (see, for example, [1], [2], [3, pp. 308–324], [4], [5]). Specifically, when a vector bundle over an $n$-dimensional manifold is associated with a representation of a structure group $G$(e.g., $GL(n)$, $O(n)$, $SO(n)$, or $Spin(n)$ in the case of spin manifolds), the covariant derivative naturally takes values in a tensor product bundle that itself admits a representation-theoretic decomposition. Extracting the $G$-irreducible components of this decomposition leads to a distinguished class of first-order operators, commonly referred to as Stein–Weiss $G$-gradients, or equivalently, generalized gradients (see [7]). These operators encode the basic $G$-irreducible components of the covariant derivative and thus provide a canonical, representation-theoretic method for constructing linear differential operators on natural tensor bundles.

Since the foundational work of Stein and Weiss (see [8]), such $G$-gradients have become central objects in global analysis. In this paper, we are primarily concerned with Stein–Weiss $O(n)$-gradients, that is, with the case $G = O(n)$. In this setting, many familiar geometric operators—such as divergences, conformal Killing operators, and others—can be realized within this framework or described as compositions of these basic gradients. For each such operator $D$, the associated second-order operator $D^*D$carries particularly rich geometric information. These operators are of interest because they often give rise to elliptic systems of differential equations with significant geometric consequences, for example in the study of Maxwell and Dirac equations.

A Weitzenböck formula expresses such an operator as the sum of the rough Laplacian (see [3])—equivalently, a second-order Stein–Weiss operator (see [5], [7])—and a curvature term determined by the representation type of the underlying bundle. These formulas not only clarify analytic properties such as ellipticity and spectral behavior, but also serve as essential tools in the study of rigidity phenomena, vanishing theorems, and curvature estimates. In particular, Weitzenböck formulas play a key role in relating local differential geometry to global topological properties via the Bochner technique (see [6]).

In dimension four, Bouguignon computed explicit Weitzenböck formulas for second-order Stein–Weiss operators acting on trace-free symmetric 2-tensors, illustrating how the special algebraic structure of curvature in four dimensions influences the analytic behavior of $D^*D$(see [3, pp. 308–324]). Branson subsequently developed general methods for the analysis of second-order Stein–Weiss operators in arbitrary dimensions and for identifying the curvature terms appearing in the corresponding Weitzenböck identities (see [1], [5], [7]).

The present paper extends these ideas to covariant symmetric trace-free tensors of arbitrary rank $p \geq 2$on a Riemannian manifold $(M, g)$. Our focus is on Stein–Weiss $O(n)$-gradients. Although the representation theory of the bundles $S^p_0(M) \coloneqq S^p_0(T^*M)$ becomes increasingly intricate as $p$grows, many of the features observed in the case $p = 2$ persist, and the higher-rank setting exhibits a systematic structure that can be made fully explicit. By analyzing the decomposition of $T^*M \otimes S^p_0(M)$ into its $O(n)$-irreducible components, we identify the natural Stein–Weiss $O(n)$-gradients acting on such tensors and derive corresponding Weitzenböck formulas for the compositions of these operators with their formal adjoints. These results generalize Bouguignon's calculations (see [3, pp. 308–324]) and provide a unified perspective on second-order operators of Stein–Weiss type in arbitrary rank.

Beyond their intrinsic geometric interest, the formulas obtained here are expected to be useful in the study of higher-rank deformation complexes, curvature-induced inequalities, and stability problems for special geometric structures. They also furnish a representation-theoretic framework in which the analytic properties of symmetric tensors can be investigated in a systematic manner.

The present article continues the authors' research on this topic, initiated in [4], [9], [10], and [11].

## 2. Stein-Weiss $O(n)$-gradients defined on symmetric traceless tensors

Let $(M, g)$ be a connected Riemannian manifold of dimension $n \geq 2$ with the Levi-Civita connection $\nabla$. Consider the bundle $S^p_0 M$ of covariant traceless symmetric $p$-tensors on $M$ for $p \geq 2$ defined by the condition tr $\varphi = \sum_{i=1,2,\dots,n} \varphi(e_i, e_i, X_3, .., X_n) = 0$ for any $\varphi \in S^p_0(T_x M)$ and $X_3, .., X_n \in T_x M$, and the orthonormal basis $e_1, \dots, e_n$ of $T_x M$ at an arbitrary point $x \in M$. The space of traceless symmetric $p$-tensors is a subbundle $S^p_0(T^*M)$ of the symmetric tensor bundle $S^p(T^*M)$. The dimension of the fiber of this bundle at any given point is:

$$dim_{\mathbb{R}} S_0^p(T_x^*M) = \binom{n+p-1}{p} - \binom{n+p-3}{p-2}$$

at each point $x \in M$.

Let Diff $\left(S_0^p M, T^*M \otimes S_0^p M\right)$ be the $C^\infty(M)$-module of first-order linear differential operators

$$D: C^\infty\left(S_0^p M\right) \to C^\infty\left(T^*M \otimes S_0^p M\right)$$

on the space $C^\infty\left(S_0^p M\right)$ of $C^\infty(M)$-sections of the bundle $S_0^p M$. For an arbitrary $p$-tensor field $\varphi \in C^\infty\left(S_0^p M\right)$, we have the following pointwise $O(n)$-irreducible decomposition of its covariant derivative (see [1]; [4]; [9] and [3, p. 329]):

$$\nabla\varphi = D_1\varphi + D_2\varphi + D_3\varphi. \tag{2.1}$$

Then in accordance with the Stein–Weiss concept (see [7] and [8]), $D_1\varphi$, $D_2\varphi$, and $D_3\varphi$ will be called Stein-Weiss $O(n)$-*gradients* or, in other words, *generalized gradients* of the symmetric tensor $\varphi \in C^\infty\left(S_0^p M\right)$ for $D_1$, $D_2$ and $D_3$ belonging to Diff $\left(S_0^p M, T^*M \otimes S_0^p M\right)$.

By [1]; [4]; [9] the first generalized gradient is the operator $D_1$: $C^\infty\left(S_0^p M\right) \to C^\infty\left(S_0^{p+1} M\right)$ which is defined by the equation:

$$D_1\varphi = \delta^*\varphi + \frac{2}{n+2(p-1)}\, g \odot \delta\varphi, \tag{2.2}$$

where

$$\left((g \odot \delta)\varphi\right)_{i_0 i_1 \dots i_p} := \frac{1}{p+1}\Big(g_{i_0 i_1}(\delta\varphi)_{i_2 \dots i_p} + g_{i_1 i_2}(\delta\varphi)_{i_3 \dots i_p i_0} + \dots +$$

$$+ g_{i_{p-1} i_p}(\delta\varphi)_{i_0 i_1 \dots i_{p-2}} + g_{i_p i_0}(\delta\varphi)_{i_1 i_2 \dots i_{p-1}}\Big);$$

$$(\delta^*\varphi)_{i_0 i_1 i_2 \dots i_{p-2} i_{p-1} i_p} := \frac{1}{p+1}\left(\nabla_{i_0}\varphi_{i_1 \dots i_p} + \nabla_{i_1}\varphi_{i_2 \dots i_p i_0} + \dots + \nabla_{i_p}\varphi_{i_0 i_1 \dots i_{p-1}}\right);$$

$$(\delta\varphi)_{i_2 \dots i_p} := -\, g^{i_0 i_1}\left(\nabla_{i_0}\varphi_{i_1}\right)_{i_2 \dots i_p},$$

where $g_{ij}$ and $\varphi_{i_1 \dots i_p}$ are local components of $g$ and an arbitrary $\varphi \in C^\infty\left(S_0^p M\right)$, respectively.

**Remark.** For $p=1$, the operator $D_1$can be identified with the conformal Killing operator, or equivalently, with the *Ahlfors operator* $S$, which maps $C^\infty(T^*M)$ to the space of trace-free symmetric tensors $C^\infty(S_0^2M)$ (see [7]). Recall that a vector field $X$ on $M$ is conformal Killing (see [12, pp. 559]) if and only if the operator $S$annihilates the one-form $\varphi$ associated with $X$ via the Riemannian metric $g$, that is, $X=\varphi^\#$. The space of all conformal Killing vector fields has dimension $n(n+1)/2+1$.

In general, if $p\geq 2$ and the tensor $\phi\in C^\infty\left(S_0^pM\right)$ lies in $\ker D_1$, then it is a traceless conformal Killing $p$-tensor (see, for example, [4], [15], and [16]). For a general $n$-dimensional Riemannian manifold with $n\geq 3$, the space of trace-free conformal Killing symmetric tensors of rank $p\geq 2$ is finite-dimensional. On a conformally flat manifold (such as Euclidean space $\mathbb{R}^n$ or the standard sphere $\mathbb{S}^n$), this dimension is explicitly known and attains the maximum possible value for any Riemannian manifold.

The dimension of $\ker D_1$ (for $p\geq 1$ and $n\geq 3$) is given by the formula (see [15]):

$$\frac{(n+p-3)!(n+p-2)!(n+2p-2)(n+2p-1)(n+2p)}{p!(p+1)!(n-2)!n!}$$

at each point $x\in M$. For a general Riemannian manifold, the dimension of this space can be lower and is bounded above by the dimension for a conformally flat manifold. As a result, we conclude that the following lemma holds.

**Lemma 1**. *Let* $(M,g)$ *be a connected Riemannian manifold of dimension* $n\geq 2$ *and* $D_1\colon C^\infty\left(S_0^pM\right)\to C^\infty\left(T^*M\otimes S_0^pM\right)$ *be the generalized gradient defined by equalities* (2.2). *Then* $\ker D_1$ *consists of the set of all traceless conformal Killing* $p$*-tensor* $(p\geq 2)$ *defined on* $(M,g)$ *and*

$$\dim_{\mathbb{R}}\ker D_1\leq\frac{(n+p-3)!(n+p-2)!(n+2p-2)(n+2p-1)(n+2p)}{p!(p+1)!(n-2)!n!}$$

*at each point* $x\in M$.

**Remark.** In particular, for the case $p=2$ it follows from (2.2) that (see also [3, p. 329])

$$(D_1\varphi)_{ijk}=\frac{1}{3}\left(\nabla_i\varphi_{jk}+\nabla_j\varphi_{ki}+\nabla_k\varphi_{ij}+\frac{2}{n+2}\left(g_{ij}(\delta\varphi)_k+g_{jk}(\delta\varphi)_i+g_{ki}(\delta\varphi)_j\right)\right).$$

In the monograph [3, p. 330], it was said that $D_1\colon C^\infty(S_0^2M) \to C^\infty(S_0^3M)$ does not have a simple geometric interpretation. On the contrary, it is known from General Relativity [12, p. 559] and [13] that a symmetric traceless 2-tensor $\varphi \in C^\infty(S_0^2M)$ is said to be a *conformal Killing tensor* in a four-dimensional space-time if it satisfies the equations

$$\nabla_i\varphi_{jk} + \nabla_j\varphi_{ki} + \nabla_k\varphi_{ij} + \frac{1}{3}\left(g_{ij}(\delta\varphi)_k + g_{jk}(\delta\varphi)_i + g_{ki}(\delta\varphi)_j\right) = 0$$

for all $i, j, k = 1, 2, 3, 4.$ Moreover, it was shown (see [14]) that the dimension of the vector space of second-order, trace-free conformal Killing tensors in a Riemannian space of dimension $n$ is bounded above by $\frac{1}{12}(n-1)(n+2)(n+3)(n+4)$ and that this is attained in flat space.

The second generalized gradient $D_2$ is defined by the equation (see [1]):

$$D_2\varphi = -\frac{n+2(p-2)}{(n+2(p-1))(n+(p-3))}\times$$
$$\times\left(g_{i_0(i_1}(\delta\varphi)_{i_2\dots i_p)} - \frac{2}{n+2(p-2)}(\delta\varphi)_{i_0(i_1 i_2\dots i_{p-2}}g_{i_{p-1}i_p)}\right) \tag{2.3}$$

Obviously, that $\ker D_2$ consists of the set of all traceless divergence-free $p$-tensors $(p \geq 2)$ defined on $(M, g)$ (see [4]). These results are consistent with the formula and conclusion given in the monograph [3, p. 329], where consider the case when $p = 2$.

Let us consider generalized gradient

$$D_2\colon C^\infty(S_0^pM) \longrightarrow C^\infty(T^*M \otimes S_0^pM),$$

such that its kernel consists of trace-free and divergence-free symmetric $p$-tensors on a Riemannian manifold $(M, g)$. We want to show that this kernel is infinite-dimensional. Namely, let $\delta$ denote the divergence operator acting on $C^\infty(S_0^pM)$ on a compact Riemannian manifold, with $(p > 1)$. This operator is underdetermined elliptic in the sense of Delay. By the results of Delay and Sunada on underdetermined elliptic operators (see [17]; [18, particularly §§7–9]; [28]), $\delta$ admits infinitely many linearly independent smooth compactly supported solutions. Therefore, the kernel of $\delta$ on trace-free symmetric tensors of rank $(p > 1)$ is infinite-dimensional. Therefore, $\ker D_2 =$ $C^\infty(S_0^p M) \cap \ker\delta$ conclude that the space of trace-free and divergence-free symmetric $p$-tensors is infinite-dimensional.

In addition, note that a traceless divergent-free conformal Killing $p$-tensor is a *Killing p-tensor* (see [4]; [12, p. 559]), which means $\delta^*\varphi = 0$ for $\varphi \in C^\infty(S_0^p M)$. Therefore, $C^\infty(S_0^p M) \cap \ker D_1 \cap \ker D_2$ consists of the set of all traceless Killing $p$-tensor ($p \geq 2$) defined on $(M, g)$. The following lemma is true.

**Lemma 2**. *Let* $(M, g)$ *be a connected Riemannian manifold of dimension* $n \geq 2$ *and* $D_2: C^\infty(S_0^p M) \to C^\infty(T^*M \otimes S_0^p M)$ *be the generalized gradient defined by equalities* (2.3). *Then* $\ker D_2$ *is an infinite-dimension vector space of the set of all traceless divergent-free p-tensors* ($p \geq 2$) *defined on* $(M, g)$. *In particular,* $\ker D_1 \cap \ker D_2$ *is a finite-dimension vector space of the set of all traceless Killing p-tensor* ($p \geq 2$) *defined on* $(M, g)$.

In particular, for $p = 2$ from (2.3) we obtain the following equations (see also [3, pp. 308-324]):

$$(D_2\varphi)_{ijk} = -\frac{n}{(n+2)(n-1)}\Big(g_{ij}\delta\varphi_k + g_{ik}\delta\varphi_j - \frac{2}{n} g_{jk}\delta\varphi_i\Big). \tag{2.4}$$

In this case, if $\varphi \in C^\infty(S_0^p M) \cap \ker D_2$, then $\delta\varphi = 0$ and it is a symmetric *transverse-traceless tensor* or, in other words, is a $TT-$ *tensor*.

**Remark.** It is well-known that $TT-$ tensors are of fundamental importance in stability analysis in General Relativity (see, for instance, [10]; [11]) and in Riemannian geometry (see [7, pp. 346-347]; [8]). As a consequence of a result of Bourguignon — Ebin — Marsden (see [7, p. 132]) the space of $TT$-tensors is an infinite-dimensional vector space on any compact (without boundary) Riemannian manifold.

**Corollary 1**. *Let* $(M, g)$ *be a connected Riemannian manifold of dimension* $n \geq 2$ *and* $D_2: C^\infty(S_0^2 M) \to C^\infty(T^*M \otimes S_0^2 M)$ *be the generalized gradient defined by equalities* (2.4). *Then* $C^\infty(S_0^2 M) \cap \ker D_2$ *is an infinite-dimensional vector space of transverse-traceless tensors* ($TT-$ *tensors*) *defined on* $(M, g)$.

The third generalized gradient $D_3$ is defined by the equation:

$$D_3\varphi = \nabla\,\varphi - D_1\varphi - D_2\varphi. \qquad (2.5)$$

In particular, for $p = 2$ we have

$$(D_3\varphi)_{ijk} = \frac{1}{3}\Big(2\left(\nabla_i\varphi_{jk} - \frac{1}{n-1}g_{jk}(\delta\varphi)_i\right) -$$

$$-\left(\nabla_j\varphi_{ki} - \frac{1}{n-1}g_{ki}(\delta\varphi)_j\right) - \left(\nabla_k\varphi_{ij} - \frac{1}{n-1}g_{ij}(\delta\varphi)_k\right)\Big)$$

If $\varphi \in \ker D_2 \cap \ker D_3$, then $\nabla\,\varphi = \delta^*\varphi$. In this case $\varphi$ is a traceless *Codazzi p-tensor* (see [4] and [19]). It is obvious that each traceless Codazzi $p$-tensor $\varphi$ has zero divergence for all $p \geq 2$. If $(M, g)$ is compact (without boundary), then $\dim_{\mathbb{R}} \ker D_2 \cap \ker D_3\varphi < \infty$, since finite dimensionality follows from the theory of elliptic operators (see [10]; [11]).

**Remark**. On a compact manifold, any overdetermined elliptic first-order operator with injective symbol has finite-dimensional kernel (see [17, p. 198]). In our case on a compact Riemannian manifold $(M, g)$, the Codazzi condition for a symmetric trace-free $(0, p)$–tensor $(p > 1)$ is given by the overdetermined elliptic first–order system $T(\varphi) = 0$, where $T(\varphi)$ is the skew–symmetrization of $\nabla\varphi$ in its first two indices. Since $T$ is overdetermined elliptic (its principal symbol is injective for all $\xi \neq 0$), standard elliptic theory implies that $\ker T$— that is, the space of Codazzi tensors of rank $p > 1$ is finite–dimensional on compact $(M, g)$.

As a result of the above, we come to the conclusion that the following lemma is valid.

**Lemma 3**. *Let $(M, g)$ be a connected Riemannian manifold of dimension $n \geq 2$. Let $D_2\colon C^\infty(S_0^2M) \to C^\infty(T^*M \otimes S_0^2M)$ and $D_3\colon C^\infty(S_0^pM) \to C^\infty(T^*M \otimes S_0^pM)$ be the generalized gradient defined by equalities* (2.4) *and* (2.5), *respectively. Then* $\ker D_2 \cap \ker D_3\varphi$ *is a finite dimensional vector space of trace-free and divergence-free Codazzi p-tensors, $p \geq 2$, defined on $(M, g)$.*

## 3. The Stein–Weiss operator $D_1^*D_1$ and the space of traceless conformal Killing tensors

Let $(M, g)$ be a connected compact (without boundary) Riemannian manifold of dimension $n \geq 2$ and $D_1\colon C^\infty(S_0^pM) \to C^\infty(S_0^{p+1}M)$ be the first-order linear differential operator defined by equalities (2.2). For the differential operator $D_1$, its formal adjoint $D_1^*$ is defined by the usual rule (see [21, p. 460])

$$\langle D_1\varphi, \Psi\rangle_{L^2} = \langle \varphi, D_1^*\Psi\rangle_{L^2}$$

for any $\varphi \in C^\infty(S_0^p M)$ and $\Psi \in C^\infty(S_0^{p+1} M)$, where $\langle\cdot,\cdot\rangle_{L^2}$ denotes the standard $L^2$-inner product

$$\langle \Phi, \Phi'\rangle_{L^2} = \int_M g(\Phi, \Phi')dv_g$$

for all smooth sections $\Phi, \Phi' \in C^\infty(\otimes^q T^*M)$ with their local components $\Phi_{i_0\dots i_q}$ and $\Phi'_{j_0\dots j_q}$. Therefore, to find $D_1^*$ we must carefully integrate by parts both terms

$$D_1\varphi = \delta^*\varphi + \frac{2}{n+2(p-1)}\, g \odot \delta\varphi.$$

We define the algebraic operator $(\operatorname{tr}\Psi)_{i_2\dots i_q} := g^{i_0 i_1}\Psi_{i_0 i_1 i_2\dots i_q}$ for an arbitrary $\Psi \in C^\infty(\otimes^{q+1} T^*M)$ with local components $\Psi_{i_0 i_1 i_2\dots i_q}$ and $(g^{ij}) = (g_{ij})^{-1}$. Then we have

$$\langle D_1\varphi, \Psi\rangle_{L^2} = \langle \delta^*\varphi, \Psi\rangle_{L^2} + \frac{2}{n+2(p-1)}\langle g \odot \delta\varphi, \Psi\rangle_{L^2} =$$

$$= \langle \varphi, \delta\,\Psi\rangle_{L^2} + \frac{2}{n+2(p-1)}\langle \delta\varphi,\ \operatorname{tr}\Psi\rangle_{L^2} = \langle \varphi, \delta\,\Psi\rangle_{L^2}$$

for any $\varphi \in C^\infty(S_0^p M)$ and $\Psi \in C^\infty(S_0^p M)$. Therefore, we can conclude that $D_1^* := \delta$. In particular, if $p = 1$, then $S^* = D_1^*$ for the conformal Killing operator $S$ (see [7]).

After this, for an arbitrary $\varphi \in C^\infty(S_0^p M)$ we substitute $D_1\varphi$ into $D_1^*$, we obtain the following:

$$D_1^* D_1\varphi = \delta\delta^*\varphi + \frac{2}{n+2(p-1)}\, \delta(g \odot \delta\varphi).$$

Using the symmetry and tracelessness of $\varphi$ and $D_1^* D_1\varphi$ we note that condition $\delta(g \odot \delta\varphi) \in C^\infty(S_0^p M)$ holds. In this case, it is directly verified that $\delta(g \odot \delta\varphi)$ can be rewritten in the form

$$\left(\delta(g\circledcirc\delta\phi)\right)_{i_1\dots i_p} = -\frac{1}{p+1}\nabla^{i_0}\Big(g_{i_0 i_1}(\delta\varphi)_{i_2\dots i_p} + g_{i_1 i_2}(\delta\varphi)_{i_3\dots i_p i_0} + \dots$$

$$+\, g_{i_{p-1} i_p}(\delta\varphi)_{i_0 i_1\dots i_{p-2}} + \, g_{i_p i_0}(\delta\varphi)_{i_1 i_2\dots i_{p-1}}\Big) =$$

$$-\frac{1}{p+1}\Big(\Big(\nabla_{i_1}(\delta\varphi)_{i_2\dots i_{p-1} i_p} + \nabla_{i_p}(\delta\varphi)_{i_1 i_2\dots i_{p-1}}\Big) -$$

$$-\Big(g_{i_1 i_2}\nabla^{i_0}\Big((\delta\varphi)_{i_3\dots i_{p-1} i_p i_0}\Big) + \dots + g_{i_{p-1} i_p}\nabla^{i_0}\Big((\delta\varphi)_{i_0 i_1 i_2\dots i_{p-2}}\Big)\Big)\Big).$$

Using that $\delta(g\circledcirc\delta\varphi)$ is trace-free, we have $g^{i_1 i_2}\left(\delta(g\circledcirc\delta\phi)\right)_{i_1 i_2 i_3\dots i_p} = 0$, and from this identity it follows that $\nabla^{i_0}(\delta\varphi)_{i_3\dots i_p i_0} = 0$. Substituting this back into the expression above eliminates all trace terms and yields

$$\left(\delta(g\circledcirc\delta\phi)\right)_{i_1\dots i_p} = -\frac{1}{p+1}\Big(\nabla_{i_1}(\delta\varphi)_{i_2\dots i_{p-1} i_p} + \nabla_{i_p}(\delta\varphi)_{i_1 i_2\dots i_{p-1}}\Big).$$

Since the left-hand side is symmetric in all indices, the right-hand side must be replaced by its full symmetrization, namely $\nabla_{(i_1}(\delta\varphi)_{i_2\dots i_p)} = (\delta^*\delta\varphi)_{i_1\dots i_p}$. Therefore,

$$\left(\delta(g\circledcirc\delta\phi)\right) = -\frac{2}{p+1}\,\delta^*\delta\varphi.$$

Combining everything, we obtain the final form:

$$D_1^* D_1\varphi = \delta\delta^*\varphi - \frac{4}{(p+1)\big(n+2(p-1)\big)}\delta^*\delta\varphi. \tag{3.1}$$

To prove ellipticity of $D_1^* D_1\varphi$, we compute its principal symbol. In normal coordinates at a point $x\in M$, for any nonzero covector $\xi\in T_x^*M$ the symbols of the divergence and its adjoint are

$$\sigma(\delta)(x,\xi)(\varphi) = \xi^{i_1}\varphi_{i_1 i_2\dots i_p}, \qquad \sigma\,(\delta^*)(x,\xi)(\varphi) = \mathrm{Sym}(\xi_{i_1}\varphi_{i_2\dots i_p}).$$

Therefore,

$$\sigma(\delta\delta^*) = \sigma(\delta)\,\sigma(\delta^*),\;\; \sigma(\delta^*\delta) = \sigma(\delta^*)\,\sigma(\delta).$$

Both compositions are *O(n)*-equivariant endomorphisms of the irreducible $O(n)$-module $S_0^p(T_xM)$. By Schur's lemma, their principal symbols are scalar:

$$\sigma(\delta\delta^*)(x,\xi) = \alpha_p\mid\xi\mid^2\;\mathrm{Id},\;\;\sigma(\delta^*\delta)(x,\xi) = \alpha_p\mid\xi\mid^2\;\mathrm{Id},$$

with $\alpha_p > 0$. Consequently,

$$\sigma(D_1^* D_1)(x,\xi) = \left(1 - \frac{4}{(p+1)(n+2(p-1))}\right) \alpha_p \mid \xi \mid^2 \ \text{Id}.$$

Since $0 < \frac{4}{(p+1)(n+2(p-1))} < 1$, the coefficient $\left(1 - \frac{4}{(p+1)(n+2(p-1))}\right) \alpha_p$ is positive, and hence $\sigma_2(D_1^* D_1)(\xi)$ is invertible for all $\xi \neq 0$. Thus $D_1^* D_1$is a *strongly elliptic* second-order operator on $S_0^p M$.

The self-adjointness and non-negativity of the elliptic operator $D_1^* D_1$ follows from the equalities (see, for example, [21, pp. 631-632]):

$$\langle D_1^* D_1 \omega, \omega' \rangle_{L^2} = \langle D_1 \omega, D_1 \omega' \rangle_{L^2} = \langle \omega, D_1^* D_1 \omega' \rangle_{L^2}, \qquad (3.2)$$

где $\langle D_1 \omega, D_1 \omega \rangle_{L^2} = \| D_1 \omega \|_{L^2}^2 \geq 0$. Moreover, the equality holds

$$C^\infty\left(S_0^p M\right) = \ker D_1^* D_1 \oplus_{L^2} \operatorname{im} D_1^* D_1.$$

As a result, we can conclude that the following theorem holds.

**Theorem 1**. *Let* $(M, g)$ *be a connected compact* (*without boundary*) *Riemannian manifold of dimension* $n \geq 2$ *and* $D_1 \colon C^\infty\left(S_0^p M\right) \to C^\infty\left(S_0^{p+1} M\right)$ *be the generalized gradient defined by equalities* (2.2). *Then the second-order Stein–Weiss operator* $D_1^* D_1 \colon C^\infty\left(S_0^p M\right) \to C^\infty\left(S_0^p M\right)$ *is given by equalities* (3.1) *and is a strongly elliptic self-adjoint and non-negative operator. Additionally, the equality holds* $C^\infty\left(S_0^p M\right) = \ker D_1^* D_1 \oplus_{L^2} \operatorname{im} D_1^* D_1$.

Since the operator $D_1^* D_1$ is strongly elliptic (see [21, p. 629]; [22, p. 383]), then its kernel $\ker D_1^* D_1$ is a finite-dimensional vector space contained in $C^\infty$ (см. [21, p. 631]). At the same time, based on the properties of elliptic operators on compact (without boundary) Riemannian manifolds (see, for example, [21, p. 632]), we conclude that the kernel of $\ker D_1^* D_1$ coincides with the kernel of the operator $D_1$, and, therefore, $\ker D_1^* D_1$ consists of the set of all traceless conformal Killing $p$-tensor ($p \geq 2$) on $(M, g)$. Hence we have

**Corollary 2**. *Let* $(M, g)$ *be a connected compact* (*without boundary*) *Riemannian manifold of dimension* $n \geq 2$ *and* $D_1 \colon C^\infty\left(S_0^p M\right) \to C^\infty\left(S_0^{p+1} M\right)$ *be the generalized gradient defined by equalities* (2.2). *Then the kernel of the second-order Stein–Weiss* *operator* $D_1^* D_1 \colon\ C^\infty\left(S_0^p M\right) \to\ C^\infty\left(S_0^p M\right)$ *is contained in* $C^\infty$ *and consists of the set of all traceless conformal Killing* $p$*-tensors* $(p\ \geq\ 2)$ *on* $(M, g)$ . *Moreover, it has finite dimension such that*

$$\dim_{\mathbb{R}} \ker D_1^* D_1 \leq \frac{(n+p-3)!(n+p-2)!(n+2p-2)(n+2p-1)(n+2p)}{p!(p+1)!(n-2)!n!}$$

*at each point* $x \in M$.

On the other hand, if $\ker D_1 = \{0\}$, then by (3.2) $\ker D_1^* D_1 = \{0\}$ and $D_1^* D_1$ becomes a strictly positive definite elliptic operator. At the same time, it is well known from [16] and [21] that on a compact Riemannian manifold $(M, g)$ of non-positive sectional curvature any trace-free conformal Killing tensor has to be parallel with respect to the Levi-Civita connection. If in addition there exists a point in $M$ where the sectional curvature of every two-plane is strictly negative, then $M$ does not carry any (non-identically zero) trace-free conformal Killing tensor (see [16]). At the same time, let us recall that the sectional curvature of a manifold is called quasi-negative if it satisfies the conditions listed above. As a result, we deduce a corollary of the above theorem.

**Corollary 3**. *Let* $(M, g)$ *be a connected compact* (*without boundary*) *Riemannian manifold of dimension* $n \geq 2$ *with quasi-negative sectional curvature and* $D_1 \colon C^\infty\left(S_0^p M\right) \to C^\infty\left(S_0^{p+1} M\right)$ *be the generalized gradient defined by equalities* (2.2). *Then the second-order Stein–Weiss operator* $D_1^* D_1 \colon\ C^\infty\left(S_0^p M\right) \to\ C^\infty\left(S_0^p M\right)$ *is a strictly positive definite elliptic operator and its kernel is zero.*

Equivalently, expression (3.1) can be rewritten using the *Sampson Laplacian*, which in this normalization is defined as (see [18]; [21, p. 356]; [23]; [24] and etc.)

$$\Delta_S := (p+1)\, \delta\delta^* - p\, \delta^*\delta.$$

Substituting this, we obtain the following form of $D_1^* D_1\, \varphi$:

$$D_1^* D_1\, \varphi = \frac{1}{p+1}\, \Delta_S\, \varphi + \frac{p}{p+1}\left(1 - \frac{2}{n+2(p-1)}\right) \delta^*\delta\, \varphi \tag{3.3}$$

This formula highlights the decomposition of $D_1^* D_1$ into the normalized Sampson Laplacian and an explicit correction term involving $\delta^* \delta \phi$, thus reflecting the geometry encoded by the Stein–Weiss gradient structure.

We conclude from (3.3) that any divergent-free traceless conformal Killing $p$-tensor (or, in other words, a traceless Killing $p$-tensor) belongs to the kernel of the Sampson Laplacian, i.e., if $\varphi \in C^\infty(S_0^p M) \cap \ker D_1 \cap \ker D_2$, then $\varphi \in \ker \Delta_S$. In the case of a compact manifold $(M, g)$, the above formula implies that

$$0 \leq \langle D_1 \varphi, D_1 \varphi \rangle_{L^2} = \frac{1}{p+1} \langle \Delta_S \varphi, \varphi \rangle_{L^2} - \frac{p(n+2(p-2))}{(p+1)(n+2(p-1))} \langle \delta\varphi, \delta\varphi \rangle_{L^2}. \tag{3.4}$$

Therefore, a traceless conformal Killing $p$-tensor belonging to the kernel of the Sampson Laplacian is divergent-free, i.e., a Killing $p$-tensor. Obviously, the converse is also true. Then we have the following corollary.

**Corollary 4**. *Let* $(M, g)$ *be a connected compact* (*without boundary*) *Riemannian manifold of dimension* $n \geq 2$ *and* $D_1 : C^\infty(S_0^p M) \to C^\infty(T^*M \otimes S_0^p M)$ *be the generalized gradient defined by equalities* (2.2). *Then the second-order Stein–Weiss operator* $D_1^* D_1$: $C^\infty(S_0^p M) \to C^\infty(S_0^p M)$ *can be rewritten in the form* (3.3) *and if* $\varphi \in Ker D_1^* D_1 \cap Ker \Delta_S$, *then* $\varphi$ *is a traceless Killing p-tensor. The converse statement is also true.*

## 4. Weitzenböck's formulas and their applications to proving vanishing theorems

At the same time, the Sampson Laplacian $\Delta_S$ acting on $C^\infty$-sections of the vector bundle $S^p M$ has the following *Weitzenböck decomposition formula* (see [24], [25], [26]):

$$\Delta_S = \nabla^* \nabla - \mathfrak{K}_p$$

where $\mathfrak{K}_p : S^p(T_x M) \to S^p(T_x M)$ is the pointwise *Weitzenböck curvature operator* which depends linearly on the Riemann curvature tensor $R$ and the Ricci tensor $Ric$ of $(M, g)$. Using the above, equations (3.3) can be rewritten in the form

$$(p+1)\, D_1^* D_1\, \varphi = \nabla^* \nabla\, \varphi - \mathfrak{K}_p(\varphi) + \frac{p(n+2(p-2))}{n+2(p-1)} \delta^* \delta\, \varphi \tag{4.1}$$

This is a *Weitzenböck formula* for the second-order Stein-Weiss operator $D_1^* D_1$ and $Q_p(\varphi, \varphi) = g\big(\mathfrak{K}_p(\varphi), \varphi\big)$ is the *Weitzenböck quadratic form* $Q_p: S_0^p(T_x M) \times S_0^p(T_x M) \to \mathbb{R}$ defined for any $\varphi \in S_0^p M$ and $x \in M$ (see, for example, [25]).

**Remark**. In particular, for $p = 1$ a similar Weitzenböck formula $S^* S = 2\delta\, d + \frac{4(n-1)}{n}\, d\, \delta - 4\, Ric$ was founded for the second-order Stein-Weiss operator or, in other word, the Ahlfors Laplacian $S^* S$, where $S := D_1$ and $d$ are the Ahlfors operator and the exterior differential defined on one-forms, respectively (see [7]). In this case, we deduce the following equation (see also [29]):

$$2\, S^* S\, \theta = \Delta_S \theta + \left(1 - \frac{2}{n}\right) d\delta\theta$$

Recall that a vector field ξ is an infinitesimal harmonic transformation in $(M, g)$ if the local one-parameter group of infinitesimal transformations generated by the vector field $\xi$ consists of harmonic transformations (see [23]). A necessary and sufficient condition for a vector field $\xi$ on a Riemannian manifold to be infinitesimal harmonic transformation in $(M, g)$ is that $\theta \in ker \Delta_S$ for $\theta^{\#} = \xi$ (see [23]). Therefore, from the above formula we conclude that any infinitesimal conformal transformation on a two-dimensional Riemannian manifold $(M, g)$ is an infinitesimal harmonic transformation (see [23]). The converse is also true. If $n \geq 3$, then from the above formula we conclude that $2\, S^* S\, \theta = \Delta_S \theta$ for an arbitrary divergence-free one-form $\theta \in C^{\infty}(T^* M)$. Moreover, any divergence-free infinitesimal harmonic transformation on an $n$-dimensional ($n \geq 3$ ) closed Riemannian manifold $(M, g)$ is an infinitesimal isometric transformation (see [23]). The converse is also true.

In conclusion of the Remark, we note that if $n \geq 3$, we can conclude from the above equation that the Ahlfors operator $S^* S: C^{\infty}(T^* M) \to C^{\infty}(T^* M)$ is not a Laplacian in the sense of the definition given in the well-known monograph [21, p. 52].

Further, (4.1) can be rewritten in integral form

$$(p + 1)\, \langle D_1 \varphi, D_1 \varphi \rangle_{L^2} = -\langle \mathfrak{K}_p(\varphi), \varphi \rangle_{L^2} + \langle \nabla \varphi, \nabla \varphi \rangle_{L^2} + \frac{p\big(n + 2(p-2)\big)}{n + 2(p-1)} \langle \delta \varphi, \delta \varphi \rangle_{L^2}. \quad (4.2)$$

If $ker\ D_1^*D_1 = \{0\}$, then from (4.1) we obtain

$$\langle \mathfrak{R}_p(\varphi), \varphi \rangle_{L^2} = \langle \nabla\varphi, \nabla\varphi \rangle_{L^2} + \frac{p(n+2(p-2))}{n+2(p-1)} \langle \delta\varphi, \delta\varphi \rangle_{L^2} \geq 0.$$

At the same time, a well-known result from [16] and [27] is that if $sec \leq 0$ everywhere on $(M, g)$, then $Q_p(\varphi, \varphi)$ is negative semidefinite for all $\varphi \in S_0^p M$ and $p \geq 2$. Using this condition, we can conclude from the above inequality that the equality $\nabla\varphi = 0$ holds for an arbitrary tensor field $\varphi \in S_0^p M$ such that $\varphi \in kerD_2 \cap kerD_3$. If in addition to the above (see [16]), there exists a point $x \in M$ where the sectional curvature of every two-plane $\pi \subset T_x M$ is strictly negative, then $Q_p(\varphi, \varphi)$ is negative definite for all $\varphi \in S_0^p(T_x M)$ and $p \geq 2$. In this case $\langle \mathfrak{R}_p(\varphi), \varphi \rangle_{L^2} < 0$. This condition contradicts the inequality written above. As a result, we conclude that the following theorem is true.

**Theorem 2**. *Let* $(M, g)$ *be a connected compact* (*without boundary*) *Riemannian manifold of dimension* $n \geq 2$ *with nonpositive sectional curvature. Let* $D_1^*D_1\colon C^\infty(S_0^p M) \to C^\infty(S_0^p M)$ *be the second-order Stein–Weiss operator for the generalized gradient* $D_1\colon C^\infty(S_0^p M) \to C^\infty(S_0^{p+1} M)$ *defined by equalities* (2.2). *Then every tensor field* $\varphi \in kerD_1^*D_1$ *is parallel, i.e.,* $\nabla\varphi = 0$. *If in addition to the above, there exists a point* $x \in M$ *where the sectional curvature of every two-plane* $\pi \subset T_x M$ *is strictly negative, then* $kerD_1^*D_1 = \{0\}$.

From (2.1) it follows that there is an identity (see also [5])

$$\nabla^*\nabla\varphi = D_1^*D_1\varphi + D_2^*D_2\varphi + D_3^*D_3\varphi.$$

Using the above identity and (4.1), we deduce the Weitzenböck formula

$$p\ D_1^*D_1\ \varphi - D_2^*D_2\varphi - D_3^*D_3\varphi = -\mathfrak{R}_p(\varphi) + \frac{p(n+2(p-2))}{n+2(p-1)} \delta^*\delta\ \varphi. \tag{4.3}$$

In this case we obtain the Weitzenböck integral formula

$$p\ \langle D_1\varphi, D_1\varphi \rangle_{L^2} - \langle D_2\varphi, D_2\varphi \rangle_{L^2} - \langle D_3\varphi, D_3\varphi \rangle_{L^2} =$$

$$= -\langle \mathfrak{R}_p(\varphi), \varphi \rangle_{L^2} + \frac{p(n+2(p-2))}{n+2(p-1)} \langle \delta\varphi, \delta\varphi \rangle_{L^2}. \tag{4.4}$$

Then using (4.4), we prove the following theorem.

**Theorem 3**. *Let* $(M, g)$ *be a connected compact* (*without boundary*) *Riemannian manifold of dimension* $n \geq 2$ *with nonnegative sectional curvature. Let* $D_2: C^\infty(S_0^p M) \to C^\infty(T^*M \otimes S_0^p M)$ *and* $D_3: C^\infty(S_0^p M) \to C^\infty(T^*M \otimes S_0^p M)$ *be the generalized gradients defined by equalities* (2.4) *and* (2.5), *respectively. Then every tensor field* $\varphi \in ker D_2^* D_2 \varphi \cap ker D_3^* D_3$ *is parallel, i.e.,* $\nabla \varphi = 0$. *In particular, every tensor field* $\varphi \in C^\infty(S_0^2 M)$ *such that* $\varphi \in ker D_3^* D_3 \cap ker D_3^* D_3$ *is zero on a connected compact Riemannian manifold* $(M, g)$ *with quasi-positive sectional curvature.*

*Proof.* Let $\varphi \in S_0^p M$ and $\varphi \in ker D_2^* D_2 \varphi \cap ker D_3^* D_3$ for $p \geq 2,$ then from (4.4) we obtain

$$- \langle \mathfrak{K}_p(\varphi), \varphi \rangle_{L^2} = p \, \langle D_1 \varphi, D_1 \varphi \rangle_{L^2} \geq 0. \qquad (4.5)$$

At the same time, a well-known result from [27] is that if $sec \geq 0$ everywhere on $(M, g)$, then $Q_p(\varphi, \varphi)$ is positive semidefinite for all $\varphi \in S_0^p M$ and $p \geq 2$. Using this condition, we can conclude from (4.5) that the equality $\nabla \varphi = 0$ holds for a tensor $\varphi \in S_0^p M$ such that field $\varphi \in ker D_2^* D_2 \varphi \cap ker D_3^* D_3$. In particular, every symmetric 2-tensor field $\varphi \in ker D_2^* D_2 \varphi \cap ker D_3^* D_3$ on a connected compact Riemannian manifold $(M, g)$ with quasi-positive sectional curvature is proportional to $g$. i.e., $\varphi = \mu \, g$ for some constant $\mu$ (see [21, p. 436]). But at present we are considering $\varphi \in S_0^2 M$ and, therefore, the equality $\varphi = 0$ holds.

A Riemannian symmetric space of compact type is an example of a compact Riemannian manifold with non-negative sectional curvature (see [21, p. 196]). Therefore, we can formulate the following corollary.

**Corollary 5.** *Let* $(M, g)$ *be a Riemannian symmetric space of compact type of dimension* $n \geq 2$. *Let* $D_2: C^\infty(S_0^p M) \to C^\infty(T^*M \otimes S_0^p M)$ *and* $D_3: C^\infty(S_0^p M) \to C^\infty(T^*M \otimes S_0^p M)$ *be the generalized gradients defined by equalities* (2.4) *and* (2.5), *respectively. Then every tensor field* $\varphi \in ker D_2^* D_2 \varphi \cap ker D_3^* D_3$ *is parallel, i.e.,* $\nabla \varphi = 0$. *In particular, if* $(M, g)$ *is a Riemannian symmetric space of compact type with quasi-positive sectional curvature, then every tensor field* $\varphi \in C^\infty(S_0^2 M)$ *such that* $\varphi \in ker D_3^* D_3 \cap ker D_3^* D_3$ *is zero.*

**Conflict of interest.** The authors declare that they have no conflict of interest.

**Funding Declaration.** The authors declare that no funds, grants, or other support were received during the preparation of this manuscript.